\documentclass{article}
\usepackage{amsmath,amsthm,amsfonts,amssymb}
\usepackage{graphicx}
\usepackage{textcomp}
\usepackage{stmaryrd}
\usepackage[cp1252]{inputenc}
\usepackage[T1]{fontenc}
\usepackage{pslatex}
\usepackage[top=2cm,bottom=2cm,right=2.5cm,left=2.5cm]{geometry}
\usepackage{wasysym}
\usepackage{amssymb}
\usepackage{tipa}
\usepackage[francais,english]{babel}
\newtheorem{lem}{Lemma}

\newtheorem{thm}{Theorem}
\newtheorem*{thmbis}{Theorem}

\newtheorem*{rem}{Remark}

\title{On a nonlinear partial integro-differential equation}

\author{Frédéric Abergel, Rémi Tachet}      
\date{November 17, 2009}

\begin{document}

\selectlanguage{english}
\large

\maketitle

\vspace{-0.5cm}
\textit{\centerline{Ecole Centrale Paris, 92295 Châtenay-Malabry, France}}
\textit{\centerline{Electronic addresses: frederic.abergel@ecp.fr, remi.tachet@ecp.fr}}
\vspace{0.5cm}

\section*{Introduction}

Financial modelling has been an area of extremely rapid growth in the past 30 years, and some extremely interesting mathematical challenges have emerged. One of the utmost importance for real-life applications  to derivatives trading is that of calibration. Similar to common situations in many areas of physics and engineering, once a model has been suggested, its parameters have to be estimated using external data. In the case of derivative modelling, those data are the liquid (tradable) options, generally known as the “vanilla” products. It is well known since the pioneering work of Litzenberger and Breeden \cite{LinBree} and its celebrated extension by Bruno Dupire \cite{dupire} that the knowledge of market data such as the prices of vanilla options across all strikes and maturities is equivalent to the knowledge of the risk-neutral marginals of the underlying stock distribution, and moreover, that there is a unique one-dimensional driftless diffusion which recovers exactly such marginals. However, it has also been well-known for almost as many years that the evolution in time of the so-called “local volatility” is not stable, thereby leading researchers and financial engineers to look for a more robust, stochastic volatility type of modelling. In this paper, we consider the calibration problem for a generic stochastic volatility model: more precisely, we address the issue of calibrating to market data a generic model with a stochastic component and a local component for the volatility process. Such models are very useful in practice, since they offer both the flexibility and realistic dynamics of stochastic volatility models, and the exact calibration properties of local volatility models. In mathematical terms, the problem we consider is a non linear partial integro-differential equation for which we are able to prove short-time existence of classical solutions under suitable assumptions.
The paper is organized as follows: Section 1 is devoted to the mathematical formulation of the problem. Section 2, to notations and statement of the main result. In Section 3, we recall some important technical results stemming from the general theory of parabolic PDE’s. Section 4 contains the proof of the main result. Finally, Section 5 is a short conclusion.

\section{The Local and Stochastic Volatility model and its calibration}

The LSV model is an extension of the Dupire local volatility model. In the simplest situation - the two-dimensional case - the dynamics of the model are given by the following system of SDE’s
\begin{align}
\frac{dS_t}{S_t} &= a(t,S_t)b(Y_t) dB^1_t + \mu_t dt  \nonumber \\
dY_t &= \alpha(t,Y_t)dB^2_t + \xi_t dt \nonumber
\end{align}
Here, $(S_t, t \geq 0)$ is the stock price process and $(Y_t,t \geq 0)$ the stochastic component of the volatility. The function $b$ simply transforms that factor into a proper volatility. $a$ is the local volatility part of the model, choosing its value properly will enable us to calibrate the vanillas of the model. $\alpha$ is the volatility of the volatility factor and $\mu$ and $\xi$ are drift terms that may depend on the state variables and on time. $B^1$ and $B^2$ are standard brownian motions with correlation $\rho$. \\
In order to fit the vanillas of this model, we write the Kolmogorov forward equation on the joint density $p(t,S,y)$ of the couple $(S_t,Y_t)$
\begin{eqnarray}
\frac{\partial p}{\partial t} - \frac{\partial^2 }{\partial S^2}(\frac{1}{2}a^2b^2S^2p) - \frac{\partial^2 }{\partial S \partial y}(\rho a b \alpha S p)  - \frac{\partial^2 }{\partial y^2}(\frac{1}{2}\alpha^2 p)) + \frac{\partial}{\partial y}(\beta p) + \frac{\partial}{\partial S}(r S p) + rp = 0  \nonumber \\
p(S,y,0)=\delta(S=S_0,y=y_0)  \nonumber
\end{eqnarray}
with $(S_0,y_0)$ the initial conditions. Taking $q = \int pdy$ the marginal of S, we get the equation
\begin{eqnarray}
\frac{\partial q}{\partial t} - \frac{\partial^2 }{\partial S^2}(\frac{1}{2}a^2S^2(\int b^2pdy)) + \frac{\partial}{\partial S}(r S q) + rq = 0  \nonumber
\end{eqnarray}
Using Dupire's results from \cite{dupire}, we know that $q$ has to solve the following equation in order to fit perfectly the vanillas of the market 
\begin{eqnarray}
\frac{\partial q}{\partial t} - \frac{\partial^2 }{\partial S^2}(\frac{1}{2}\sigma_D^2S^2 q) + \frac{\partial}{\partial S}(r S q) + rq = 0  \nonumber
\end{eqnarray}
where $\sigma_D$ is Dupire's local volatility and contains the information about the vanillas we want to reproduce.
We identify the terms in this last formula. This gives us the value of $ a^2(t,S) = \sigma_D^2(t,S) \frac{q}{\int b^2 pdy} = \sigma_D^2(t,S) \frac{\int pdy}{\int b^2 pdy}$.
Eventually, the joint density that calibrates the smile of our model is solution of the nonlinear partial integro-differential equation
\begin{eqnarray}
\frac{\partial p}{\partial t} - \frac{\partial^2 }{\partial S^2}(\frac{1}{2}\sigma_D^2 b^2S^2 \frac{\int pdy}{\int b^2 pdy} p) - \frac{\partial^2 }{\partial S \partial y}(\rho \sigma_D b \alpha S  (\frac{\int pdy}{\int b^2 pdy})^{\frac{1}{2}} p)  - \frac{\partial^2 }{\partial y^2}(\frac{1}{2}\alpha^2 p)) \nonumber \\
+ \frac{\partial}{\partial y}(\beta p) + \frac{\partial}{\partial S}(r S p) + rp = 0  \nonumber 
\end{eqnarray}
The rest of this paper is devoted to the study of a more general n-dimensional version of this equation.

\section{Generalized equation and notations}

Throughout this article, we denote by $0 < t \leq T$ the time-variable and by $x=(x_1,x_2,...,x_n) \in \Omega \subset \mathbb{R}^n$ the n-dimensional space variable where $\Omega$ is an open subset with a sufficiently smooth boundary (we will precise this notion later). When we consider the equation from a financial point of view, the first-variable $x_1$ stands for the spot and the last $n-1$ for the volatility. Hence, we write $S=x_1$, $y=(x_2,...,x_n)$. We also denote by $D^T = ]0,T[$ \texttimes \hspace{.03cm} $\Omega$ the domain of definition and by $B = \{0\}$ \texttimes \hspace{.03cm} $\Omega$, $B^T = \{T\}$ \texttimes \hspace{.03cm} $\Omega$ and $C^T = ]0,T[$ \texttimes \hspace{.03cm} $\partial \Omega$ the different parts of the boundary. Given the particular part played by the spot, we consider $\Omega_S = \{S \in \mathbb{R} / \exists y \in \mathbb{R}^{n-1}, (S,y) \in \Omega\}$ and $\forall S \in \Omega_S, \Omega^S_y = \{y \in \mathbb{R}^{n-1} / (S,y) \in \Omega\}$.  
We are interested in the following equation:
\begin{eqnarray}
O(p) := \frac{\partial p}{\partial t} - \frac{\partial^2 }{\partial S^2}(\rho_{11}\alpha_1^2 I(p) p) - \displaystyle {\sum_{i=2}^{n}} \frac{\partial^2 }{\partial S \partial x_i}(\rho_{1i} \alpha_1 \alpha_i  \sqrt{I(p)} p) \nonumber && \\  - \displaystyle {\sum_{i,j=2}^{n}} \frac{\partial^2 }{\partial x_i \partial x_j}(\rho_{ij} \alpha_i \alpha_j p)  + \displaystyle {\sum_{i=1}^{n}} \frac{\partial}{\partial x_i}(\beta_i p) + \gamma p = 0 &on& D^T \cup B^T
\label{eq_gen_div}
\end{eqnarray} 
where $(\rho_{ij})_{1 \leq i,j \leq n}$ is a correlation matrix ie is symetric positive definite and verifies: $\rho_{ii} = \frac{1}{2}$ for all $i$ and $-\frac{1}{2} < \rho_{ij} < \frac{1}{2}$ for $i\neq j$. We add the boundary condition $p =  \Psi$ on $\overline{B} \cup C^T$ with $\Psi$ constant on $C^T$. We also let $p_0$ denote the function $p_0 (t,S,y) = \Psi (S,y)$.
The complexity of this equation stems from the following integral term:
\begin{eqnarray}
I(p)(t,S) = \frac{\int_{\Omega^S_y}{p(t,S,x_2,...,x_n)dx_2...dx_n}}
{\int_{\Omega^S_y}{b^2(x_2,...,x_n)p(t,S,x_2,...,x_n)dx_2...dx_n}}
= \frac{\int_{\Omega^S_y}{p(t,S,y)dy}}{\int_{\Omega^S_y}{b^2(y)p(t,S,y)dy}} 
\label{terme_int}
\end{eqnarray}

\noindent
(\ref{eq_gen_div}) belongs to the class of nonlinear, parabolic and nonlocal equations. An interesting reference concerning that kind of equations is \cite{alib}. However, our case doesn't fall under the scope of that paper: the operator $I(.)$ is not defined on $C_b(\overline{D^T})$.

\noindent
Let us now make a few remarks about our particular equation:
\begin{enumerate}
\item in the case of an equation with no $I(p)$ term, it becomes a classic linear equation of parabolic type. That kind of equation has been properly solved for quite some time now, see \cite{friedman} or \cite{lady}.
\item when $b$ is constant, the problem is reduced to the previous remark. This observation is the key to our resolution method. First, we suppose that b does not vary too much and approximate the nonlocal term $I(p)$ with a suitable constant. We then isolate the error made during this process in the second term and use a fixed point method to solve the new equation. 
\item in order to use this method and the results from \cite{friedman}, one has to assume that the coefficients of the equation belong to Hölder spaces $H^{k,h,h/2}$ (we shall define them in the preliminaries). 
\item the question whether I is properly defined is natural. To answer it, we have to prove that $\int_{\Omega^S_y}{b^2(y)p(t,S,y)dy}$ is bounded away from $0$. To do it, we assume that b is non-negative and that $\Psi$ the initial condition is strictly positive. By restricting ourselves to short times, we are sure that p is not too far from its initial condition and thus is strictly positive. 
\item from a financial viewpoint, it is natural to consider a domain $\Omega$ cylindrical with respect to the spot. However, since it may become very challenging to study a PDE on a domain with corners, we reduce our study to domains with a $S$-section depending on $S$.
\end{enumerate} 
\noindent
The theorem we will prove requires the following assumptions on $b$ and $\Psi$. 
\begin{description}
\item (H1) $b \in C^{1}(\mathbb{R}^{n-1})$, $\exists (\delta_1,\delta_2) \in \mathbb{R}^2$, $0 < \delta_1 \leq b \leq \delta_2$ on $\mathbb{R}^{n-1}$
\item (H2) $\forall 2 \leq i \leq n, |\frac{\partial b^2}{\partial x_i}| \leq b^*$ on $\Omega_y$ where $b^*$ is a constant we will choose later
\item (H3) $\Psi$ is strictly positive and in $H^{2,h,h/2}$. This gives us two results on $p_0$. First, $p_0$ belongs to $H^{2,h,h/2}(D^T)$ and second
\begin{eqnarray}
0 \hspace{3pt} < \hspace{3pt} \underline{p_0} = inf \hspace{3pt} p_0 \hspace{3pt} \leq \hspace{3pt} sup \hspace{3pt} p_0 = \overline{p_0} 
\nonumber
\end{eqnarray}
\item (H4) $O(\Psi) = 0$ on $\partial B$ in a sense described in the preliminaries
\end{description}

\noindent 
Under the previous assumptions, we have the following result:
\begin{thm}\label{theo_gen_div}
If the $\alpha_i$ belong to $H^{2,h,h/2}(D^T)$, are positive and bounded away from 0 by a stricly positive constant e, if the $\beta_i$ are in $H^{1,h,h/2}(D^T)$ and if $\gamma$ belongs to $H^{0,h,h/2}(D^T)$, then, for $b^*$ small enough, there exists $0 < T^* \leq T$ and a solution of the equation (\ref{eq_gen_div}) on $D^{T^*} \cup B^{T^*}$.
\end{thm}
\noindent
The rest of the paper will be devoted to the proof of Theorem 1 above.

\section{Preliminaries}

In this section, we write $\partial_{x_i}$ and $\frac{\partial}{\partial x_i}$ without distinction. Let us start as in \cite{friedman} and \cite{lady} with the following notion of distance $d(P,Q) = [|x - x'|^2 + |t - t'|]^{1/2}$ where $P = (t,x)$ and $Q = (t',x')$ belong to $D^T$ and $|x|$ is the norm of the n-dimensionnal vector $x$. Given such a metric d, we can define the concept of Hölder continuity. For a function u, we write:
\begin{eqnarray}
|u|^{D^T}_0 = \underset{D^T}{sup} |u| \hspace{1.5cm} H^{D^T}_h(u) = \underset{P,Q \in D^T}{sup} \frac{|u(P) - u(Q)|}{d(P,Q)^h} \hspace{1.5cm} |u|^{D^T}_h = |u|^{D^T}_0 + H^{D^T}_h(u) \nonumber
\end{eqnarray}

\noindent
$H^{D^T}_h(u) < \infty$ if and only if $u$ is uniformly hölder (exponent $h$) in ${D^T}$. We denote by $H^{0,h,h/2}(D^T)$ the set of all functions u for which $|u|^{D^T}_h < \infty$. Now, if all the derivatives used in the equation exist, we write for $k \leq 2$:
\begin{eqnarray}
|u|^{D^T}_{k+h} = |u|^{D^T}_{h} + \Sigma |\partial_{x} u|^{D^T}_{h} + ...  + \Sigma |\partial^k_{x} u|^{D^T}_{h} + |\partial_{t} u|^{D^T}_{h} \label{norm_hol}
\end{eqnarray}
where the sums are taken over all the partial derivatives of the indicated order. We denote by $H^{k,h,h/2}(D^T)$ the set of all functions u for which $|u|^{D^T}_{k+h} < \infty$. It is a Banach space and an algebra with the norm given by definition \ref{norm_hol}. Indeed, for all u,v in $H^{k,h,h/2}(D^T)$, we have:
\begin{eqnarray}
|uv|^{D^T}_{k+h} \leq |u|^{D^T}_{k+h} |v|^{D^T}_{k+h}
\end{eqnarray}
We can now make the assumptions about $D^T$ more precise: for every point $Q$ of $\overline{C^T}$, there exists an (n+1)-dimensional neighborhood V such that $V \cap \overline{C^T}$ can be represented, for some i ($1 \leq i \leq n$), in the form $x_i = r(t,x_1,...,x_{i-1},x_{i+1},...,x_n)$ with $r$, $\partial_x r$, $\partial^2_{xx} r$, $\partial_t r$ Hölder continuous (exponent h) and $\partial^2_{xt} r$, $\partial^2_{tt} r$ simply continous. \newline 
We also have to consider functions $\psi$ defined on $\overline{B} \cup C^T$. Such a function $\psi$ is said to belong to $H^{k,h,h/2}$ if there exists a $\Psi$ in $H^{k,h,h/2}(D^T)$ such that $\Psi = \psi$ on $\overline{B} \cup C^T$. We then define $|\psi|_{k+h} =$ \vspace{.2cm}inf $|\Psi|_{k+h}^{D^T}$ where the inf is taken with respect to all the $\Psi$'s in $H^{k,h,h/2}(D^T)$ which coincide with $\psi$ on $\overline{B} \cup C^T$. This process defines a norm on $H^{k,h,h/2}$.

\noindent
The following results will be useful in the proof of our result. They concern the PDE:
\begin{eqnarray}
\label{eq_frie}
Lu := \frac{\partial u}{\partial t} - \displaystyle {\sum_{i,j=1}^{n}} a_{ij}(x,t) \frac{\partial^2 u}{\partial x_i \partial x_j}  + \displaystyle {\sum_{i=1}^{n}} b_i(x,t) \frac{\partial u}{\partial x_i} + c(x,t) u = f(x,t) & on & D^T \cup B^T \\
u = \psi & on & \overline{B} \cup C^T \nonumber
\end{eqnarray} 

\noindent
We shall need the assumptions:
\begin{itemize}
\item the coefficients of the operator L belong to $H^{0,h,h/2}(D^T)$, let $K_1$ be a bound on their norm
\item for all $(x,t)$ in $D^T$ and for all $\xi \in \mathbb{R}^n$, $\displaystyle{ \sum^n_{i,j=1}} a_{ij}(x,t) \xi_i \xi_j \geq K_2 \mid \xi \mid^2$ $(K_2 > 0)$ 
\item $\psi \in H^{2,h,h/2}$ and $|f|_h^{D^T} < \infty $ 
\end{itemize}
\noindent
In addition, given the assumption about $D^T$, if we consider a function $\psi \in H^{2,h,h/2}$, for any extension $\Psi$ of $\psi$, $\partial_t \Psi$ is uniquely defined (by continuity) on the boundary $\partial B$ of $B$, and the definition is independent of $\Psi$. We denote this function (on $\partial B$) by $\partial_t \psi$. The other terms of $L \psi$ are also uniquely defined (by continuity) on $\partial B$. Thus, the quantity $L \psi$ is well-defined on $\partial B$.

\begin{thm}
\label{res_frie}
Under the previous assumptions and if $L\psi = f$ on $\partial B$, there exists a unique solution of the equation \ref{eq_frie}, this solution belongs to $H^{2,h,h/2}(D^T)$ and we have the Schauder inequality (with $K_{H^2}$ depending only on $K_1$, on $K_2$, on $h$ and on $D^T$)
\begin{eqnarray}
|u|_{2+h}^{D^t} &\leq& K_{H^2} (|\psi|_{2+h} + |f|_h^{D^t}) \label{schau_2} 
\end{eqnarray}
Furthermore, if $\psi = 0$, we can write a bound containing the time on the supremum of the solution
\begin{eqnarray}
|u|_{0}^{D^t} &\leq& t K_{0} |f|_{0}^{D^t}  \label{schau_inf} 
\end{eqnarray}
where $K_{0}$ only depends on $K_1$, on $K_2$, on $h$ and on $\Omega$.
\end{thm}

Proof: the first part of the result is classic, its proof can be found in \cite{friedman}. As to the result with $\psi = 0$, which is more original, one needs a result from \cite{lady} about volume potentials and representation of solutions of parabolic equations. It is the theorem (16.2) of section IV.16 we shall use. One reads that the solution of the equation \ref{eq_frie} with $\psi = 0$ can be written as 
\begin{eqnarray}
u(x,t) = \int^t_0 d\tau \int_\Omega G(x,z,t,\tau) f(z, \tau) dz \nonumber
\end{eqnarray}
where $G$ is the Green's function for the operator L and verifies
\begin{eqnarray}
|G(x,y,t,\tau)| \leq K (t - \tau)^{-\frac{n}{2}} exp(-K' \frac{|x-y|^2}{t - \tau})
\end{eqnarray}
with $K$ and $K'$ two constants depending on the data of the problem. Using both these results, we get, for all $t' \leq t$ and $x \in \Omega$
\begin{eqnarray}
|u(x,t')| \leq t |f|_{0}^{D^t} \int^t_0 d\tau \int_\Omega K (t - \tau)^{-\frac{n}{2} - 1} exp(-K' \frac{|x-y|^2}{t - \tau}) dz \leq t |f|_{0}^{D^t} K_{0} \nonumber
\end{eqnarray} 
where $K_{0}$ depends on $K_1$, on $K_2$, on $h$ and on $D^T$. 

\section{Proof of Theorem \ref{theo_gen_div}}

We are interested in the equation
\begin{eqnarray}
\frac{\partial p}{\partial t} - \frac{\partial^2 }{\partial S^2}(\rho_{11}\alpha_1^2 I(p) p) - \displaystyle {\sum_{i=2}^{n}} \frac{\partial^2 }{\partial S \partial x_i}(\rho_{1i} \alpha_1 \alpha_i  \sqrt{I(p)} p) \nonumber && \\  - \displaystyle {\sum_{i,j=2}^{n}} \frac{\partial^2 }{\partial x_i \partial x_j}(\rho_{ij} \alpha_i \alpha_j p)  + \displaystyle {\sum_{i=1}^{n}} \frac{\partial}{\partial x_i}(\beta_i p) + \gamma p = 0 &on& D^T \cup B^T \nonumber
\end{eqnarray} 
and want to prove the

\begin{thmbis}
If the $\alpha_i$ belong to $H^{2,h,h/2}(D^T)$, are positive and bounded away from 0 by a stricly positive constant e, if the $\beta_i$ are in $H^{1,h,h/2}(D^T)$ and if $\gamma$ belongs to $H^{0,h,h/2}(D^T)$, then, for $b^*$ small enough, there exists $0 < T^* \leq T$ and a solution of the equation (\ref{eq_gen_div}) on $D^{T^*} \cup B^{T^*}$.
\end{thmbis}

Proof:
\newline
The assumption (H2) gives us some control over the variations of b. Let us denote by $\underline {b} = b(y_0)$ a strictly positive value taken by b (with $y_0 \in \Omega^S_y$ for some arbitrary $S \in \Omega_S$). We use the assumption on b to approximate the integral term $I(p)$ with $1/{\underline {b}}^2$, the gap between those two quantities is quantified with the
\begin{lem}\label{lemme} 
There exists a constant $K_{b}$ (depending only on $h$, $n$, $\delta_1$, $\delta_2$, $\underline{p_0}$ and $\Omega$) and a polynomial function P strictly positive and increasing on $\mathbb{R}^*_+$ such that $\forall p \in H^{2,h,h/2}(D^T)$ verifiying $\underline{p_0} \leq p$, we have
\begin{eqnarray}
|I(p) - \frac{1}{{\underline {b}}^2}|^{D^t}_{2+h} + |\sqrt{I(p)} - \frac{1}{\underline{b}}|^{D^t}_{2+h}
 \leq b^* K_{b} P(|p|^{D^t}_{2+h}). \nonumber
\end{eqnarray}
\end{lem}

\begin{rem}
As a consequence of this lemma, we see that $\forall p \in H^{2,h,h/2}(D^T)$ verifiying $\underline{p_0} \leq p$, $I(p)$ belongs to $H^{2,h,h/2}(D^T)$.
\end{rem}

\noindent
We then write the equation as
\begin{eqnarray}
\frac{\partial p}{\partial t} - \frac{\partial^2 }{\partial S^2}(\rho_{11}\alpha_1^2 \frac{1}{{\underline {b}}^2} p) - \displaystyle {\sum_{i=2}^{n}} \frac{\partial^2 }{\partial S \partial x_i}(\rho_{1i} \alpha_1 \alpha_i  \frac{1}{\underline {b}} p)  - \displaystyle {\sum_{i,j=2}^{n}} \frac{\partial^2 }{\partial x_i \partial x_j}(\rho_{ij} \alpha_i \alpha_j p)  + \displaystyle {\sum_{i=1}^{n}} \frac{\partial}{\partial x_i}(\beta_i p) + \gamma p \nonumber && \\ = \frac{\partial^2 }{\partial S^2}(\rho_{11}\alpha_1^2 (I(p) - \frac{1}{{\underline {b}}^2}) p) + \displaystyle {\sum_{i=2}^{n}} \frac{\partial^2 }{\partial S \partial x_i}(\rho_{1i} \alpha_1 \alpha_i  (\sqrt{I(p)} - \frac{1}{\underline {b}}) p) && \nonumber
\end{eqnarray}
To solve this equation, we apply a fixed point method and use the lemma \ref{lemme} to get an upper bound on the second term. \\
We take a real number $x \geq |p_0|_{2+h}^{D^T}$ and $t \in \mathbb{R}_+^*$ and let $X^t_x$ denote the set
\begin{equation}
X^t_x = \{p \in H^{2,h,h/2}(D^t), \hspace{3pt} |p|_{2+h}^{D^t} \leq x, \hspace{3pt} \frac{\underline{p_0}}{2} \leq p \leq \overline{p_0} + \frac{\underline{p_0}}{2}, \hspace{3pt} p = \Psi \hspace{4pt} on \hspace{4pt} \overline{B} \cup C^T \} \nonumber
\end{equation}
The set $X^t_x$ clearly contains the function $p_0$. We then consider the application M which takes a function $u \in X$ and sends it on $v \in H^{2,h,h/2}(D^T)$ solution of the equation
\begin{eqnarray}
O'v := \frac{\partial v}{\partial t} - \frac{\partial^2 }{\partial S^2}(\rho_{11}\alpha_1^2 \frac{1}{{\underline{b}}^2} v) - \displaystyle {\sum_{i=2}^{n}} \frac{\partial^2 }{\partial S \partial x_i}(\rho_{1i} \alpha_1 \alpha_i  \frac{1}{\underline{b}} v)  - \displaystyle {\sum_{i,j=2}^{n}} \frac{\partial^2 }{\partial x_i \partial x_j}(\rho_{ij} \alpha_i \alpha_j v)  + \displaystyle {\sum_{i=1}^{n}} \frac{\partial}{\partial x_i}(\beta_i v) + \gamma v \nonumber && \\ = \frac{\partial^2 }{\partial S^2}(\rho_{11}\alpha_1^2 (I(u) - \frac{1}{{\underline{b}}^2}) u) + \displaystyle {\sum_{i=2}^{n}} \frac{\partial^2 }{\partial S \partial x_i}(\rho_{1i} \alpha_1 \alpha_i  (\sqrt{I(u)} - \frac{1}{\underline{b}}) u) && 
\label{point_fixe_div_1}
\end{eqnarray}
with the boundary condition $v = \Psi$ on $\overline{B} \cup C^T$. The existence of v is given by Theorem \ref{res_frie}. Indeed, the coefficients of this equation belong to the appropriate spaces and because of (H4) the necessary condition
\begin{eqnarray}
O' \psi =  \frac{\partial^2 }{\partial S^2}(\rho_{11}\alpha_1^2 (I(\Psi) - \frac{1}{\underline{b}^2}) \Psi) + \displaystyle {\sum_{i=2}^{n}} \frac{\partial^2 }{\partial S \partial x_i}(\rho_{1i} \alpha_1 \alpha_i  (\sqrt{I(\Psi)} - \frac{1}{\underline{b}}) \Psi) \nonumber
\end{eqnarray} 
on $\partial B$ is verified. It remains to prove that this operator is elliptic: let $(\xi_i)_{1 \leq i \leq n}$ be n real numbers and $(t,S,y) \in D^T$. We write $f_1 = \frac{\alpha_1}{\underline{b}}$ and $f_i = \alpha_i$ for $i \geq 2$, we have
\begin{eqnarray}
\displaystyle {\sum_{i,j=1}^n} \rho_{ij} f_i (t,S,y) f_j (t,S,y) \xi_i \xi_j
&\geq& K_{\rho} \displaystyle {\sum_{i=1}^n} f_i^2(t,S,y) \xi_i^2 \geq K_{\rho} e^2 \displaystyle {\sum_{i=1}^n} \xi_i^2 \nonumber 
\end{eqnarray}
where the existence of $K_{\rho}$ is a consequence of $\rho$ being a positive definite matrix.
This proves the ellipticity of the operator, v exists and belongs to $H^{2,h,h/2}(D^T)$. We now want to show that for suitable $x$ and $t$, v belongs to $X^t_x$ ie that
\begin{eqnarray}
|v|_{2+h}^{D^t} \leq x  \hspace{2cm} \frac{\underline{p_0}}{2} \leq v \leq \overline{p_0} + \frac{\underline{p_0}}{2} \nonumber
\end{eqnarray}
For the first inequality, we apply \ref{schau_2}


\begin{eqnarray}
|v|_{2+h}^{D^t}
&\leq& K_{H^2} (|\psi|_{2+h}^{D^t} + |\frac{\partial^2 }{\partial S^2}(\rho_{11}\alpha_1^2 (I(u) - \frac{1}{{\underline{b}}^2}) u) + \displaystyle {\sum_{i=2}^{n}} \frac{\partial^2 }{\partial S \partial x_i}(\rho_{1i} \alpha_1 \alpha_i  (\sqrt{I(u)} - \frac{1}{\underline{b}}) u)|_h^{D^t}) \nonumber \\
&\leq& K_{H^2} (|\psi|_{2+h}^{D^T} + |\rho_{11}\alpha_1^2 (I(u) - \frac{1}{{\underline{b}}^2}) u |_{2+h}^{D^t}+ \displaystyle {\sum_{i=2}^{n}} |\rho_{1i} \alpha_1 \alpha_i  (\sqrt{I(u)} - \frac{1}{\underline{b}}) u|_{2+h}^{D^t}) \nonumber \\
&\leq& K_{H^2} (|\psi|_{2+h}^{D^T} + (|\rho_{11}\alpha_1^2|_{2+h}^{D^T} + \displaystyle {\sum_{i=2}^{n}} |\rho_{1i} \alpha_1 \alpha_i|_{2+h}^{D^T}) (|I(u) - \frac{1}{{\underline{b}}^2}|_{2+h}^{D^t} + |\sqrt{I(u)} - \frac{1}{\underline{b}}|_{2+h}^{D^t}) |u|_{2+h}^{D^t}) \nonumber \\
&\leq& K_{H^2} (|\psi|_{2+h}^{D^T} + b^* K' P(|u|_{2+h}^{D^t}) |u|_{2+h}^{D^t}) \label{in_1} 
\end{eqnarray}

\noindent
where $K' = (|\rho_{11}\alpha_1^2|_{2+h}^{D^T} + \displaystyle {\sum_{i=2}^{n}} |\rho_{1i} \alpha_1 \alpha_i|_{2+h}^{D^T}) K_b$ (we apply lemma \ref{lemme} for the last line). \\

\noindent
We remember that u belongs to $X^t_x$, thus $|u|_{2+h}^{D^t} \leq x$ and then
\begin{eqnarray}
|v|_{2+h}^{D^t} \leq K_{H^2} (|\psi|_{2+h}^{D^T} + b^* K' P(x) x) \nonumber 
\end{eqnarray}
Taking $x^* = max (K_{H^2} (|\psi|_{2+h}^{D^T} + 1), |p_0|_{2+h}^{D^T}$) and $b^* \leq \frac{1}{K' P(x^*) x^*}$, \ref{in_1} gives us $|v|_{2+h}^{D^T} \leq x^*$. \\
It remains to prove that $\frac{\underline{p_0}}{2} \leq v \leq \overline{p_0} + \frac{\underline{p_0}}{2}$.
Let us write $\tilde{v} = v - p_0$. It is clear that $\tilde{v}$ verifies 
\begin{eqnarray}
O' \tilde{v} = O' p_0 + \frac{\partial^2 }{\partial S^2}(\rho_{11}\alpha_1^2 (I(u) - \frac{1}{{\underline{b}}^2}) u) + \displaystyle {\sum_{i=2}^{n}} \frac{\partial^2 }{\partial S \partial x_i}(\rho_{1i} \alpha_1 \alpha_i  (\sqrt{I(u)} - \frac{1}{\underline{b}}) u) \nonumber
\end{eqnarray} 
on $D^T \cup B^T$ with $\tilde{v} = 0$ on $\overline{B} \cup C^T$ (here we use the fact that $\Psi$ is constant on $C^T$). \\
We now apply the second part of Theorem \ref{res_frie}, the inequality \ref{schau_inf}, to this function $\tilde{v}$
\begin{eqnarray}
|\tilde{v}|_{0}^{D^t}
&\leq& t K_{0} |O' p_0 + \frac{\partial^2 }{\partial S^2}(\rho_{11}\alpha_1^2 (I(u) - \frac{1}{{\underline{b}}^2}) u) + \displaystyle {\sum_{i=2}^{n}} \frac{\partial^2 }{\partial S \partial x_i}(\rho_{1i} \alpha_1 \alpha_i  (\sqrt{I(u)} - \frac{1}{\underline{b}}) u)|_{0}^{D^t} \nonumber \\
&\leq& t K_{0} (|O' p_0|_{0}^{D^T} + |\frac{\partial^2 }{\partial S^2}(\rho_{11}\alpha_1^2 (I(u) - \frac{1}{{\underline{b}}^2}) u) + \displaystyle {\sum_{i=2}^{n}} \frac{\partial^2 }{\partial S \partial x_i}(\rho_{1i} \alpha_1 \alpha_i  (\sqrt{I(u)} - \frac{1}{\underline{b}}) u)|_{h}^{D^t}) \nonumber \\
&\leq& t K_{0} (|O' p_0|_{0}^{D^T} + 1) \label{in_2}
\end{eqnarray}
Taking $T^* K_{0} (|O' p_0|_{0}^{D^T} + 1) = \frac{\underline{p_0}}{2}$, we get $|\tilde{v}|_{0}^{D^t} \leq \frac{\underline{p_0}}{2}$. Eventually, since $v = p_0 + \tilde{v}$, the last inequality is proved and v belongs to $X^{T^*}_{x^*}$. The application M maps $X^{T^*}_{x^*}$ into itself. \\
Using this statement, we construct a bounded sequence $(p_n)_{n \in \mathbb{N}}$ of functions belonging to $X^{T^*}_{x^*}$
\begin{itemize}
\item $p_0$ has been previously defined
\item by induction, we write $p_{n+1} = M(p_n)$
\end{itemize}
By construction, we have $\forall n \in \mathbb{N}, |p_n|_{2+h}^{D^T} \leq x^*$. Repeated applications of the Ascoli-Arzelà theorem give us a function $p \in C^2(D^{T^*})$ limit in $C^2(D^{T^*})$ of a subsequence of $p_n$. Since
\begin{eqnarray}
sup \{\frac{\mid \frac{\partial^2 p_n}{\partial x_i \partial x_j}(x,t)- \frac{\partial^2 p_n}{\partial x_i \partial x_j} (x',t') \mid}{(\mid x - x' \mid^2 + \mid t - t' \mid)^{\alpha / 2}}; (x,t),(x',t') \in D^{T^*}\} \leq x \nonumber
\end{eqnarray}
We have
\begin{eqnarray}
sup \{\frac{\mid \frac{\partial^2 p}{\partial x_i \partial x_j}(x,t)- \frac{\partial^2 p}{\partial x_i \partial x_j} (x',t') \mid}{(\mid x - x' \mid^2 + \mid t - t' \mid)^{\alpha / 2}}, (x,t);(x',t') \in D^{T^*}\} \leq x \nonumber
\end{eqnarray}
And this computation being true for all the derivatives appearing in the norm $H^{2,h,h/2}$, we find that $p \in H^{2,h,h/2}$. \\
The last step of the proof is to take the limit in \ref{point_fixe_div_1}. The only result needed is $I(p_n) \rightarrow I(p)$. Since $p_n \in X^{T^*}_{x^*}$, the denominator is bounded away from 0. Two applications of the dominated convergence theorem give us the convergence we need. Thus, it is possible to take the limit in \ref{point_fixe_div_1} which gives us

\begin{eqnarray}
\frac{\partial p}{\partial t} - \frac{\partial^2 }{\partial S^2}(\rho_{11}\alpha_1^2 \frac{1}{{\underline{b}}^2} p) - \displaystyle {\sum_{i=2}^{n}} \frac{\partial^2 }{\partial S \partial x_i}(\rho_{1i} \alpha_1 \alpha_i  \frac{1}{\underline{b}} p)  - \displaystyle {\sum_{i,j=2}^{n}} \frac{\partial^2 }{\partial x_i \partial x_j}(\rho_{ij} \alpha_i \alpha_j p)  + \displaystyle {\sum_{i=1}^{n}} \frac{\partial}{\partial x_i}(\beta_i p) + \gamma p \nonumber && \\ = \frac{\partial^2 }{\partial S^2}(\rho_{11}\alpha_1^2 (I(p) - \frac{1}{{\underline{b}}^2}) p) + \displaystyle {\sum_{i=2}^{n}} \frac{\partial^2 }{\partial S \partial x_i}(\rho_{1i} \alpha_1 \alpha_i  (\sqrt{I(p)} - \frac{1}{\underline{b}}) p) && \nonumber 
\end{eqnarray}
That concludes the proof of the theorem, p is solution of our equation. 

\noindent
Let us now prove lemma \ref{lemme}. By definition, we have:
\begin{eqnarray}
|I(p) - \frac{1}{{\underline {b}}^2}|^{D^T}_{2+h} = |\frac{\int{pdy}}{\int{b^2pdy}} - \frac{1}{{\underline {b}}^2}|^{D^T}_{2+h} &=& |\frac{\int{p(\underline {b}^2 - b^2)dy}}{\underline {b}^2 \int{b^2pdy}}|^{D^T}_{2+h}
\nonumber \\
&\leq& \frac{1}{\underline {b}^2} |\frac{1}{\int{b^2pdy}}|^{D^T}_{2+h} |\int{p(\underline {b}^2 - b^2)dy}|^{D^T}_{2+h} \nonumber
\end{eqnarray}
Let us compute one after another the terms appearing in this norm (we remember that those functions only depend on $t$ and $S$). Let $(t,S)$ belong to {$]0,T[$ \texttimes \hspace{.02cm} $\Omega_S$}. We have
\begin{eqnarray}
|\int_{\Omega_y^S}{p(t,S,y) (\underline {b}^2 - b^2(y))dy}| && \leq \int_{\Omega_y^S}{p(t,S,y) |\underline {b}^2 - b^2(y)| dy}  \leq |p|^{D^T}_{0} \int_{\Omega_y^S}{|\underline {b}^2 - b^2(y)| dy} \nonumber \\
&& \leq |p|^{D^T}_{0} \int_{\Omega_y^S}{|b^2(y_0) - b^2(y)| dy} \leq b^* (n-2) \int_{\Omega_y^S}{|y_0 - y| dy} |p|^{D^T}_{0} \nonumber \\
&& \leq b^* K |p|^{D^T}_{0} \nonumber
\end{eqnarray}
here and in the rest of the proof, K stands for some constant depending only on the data of the problem ($\Omega$, $n$, $\delta_1$...) but not on $p$ nor on $b^*$.
We get the last line from the following computation where $y = (x_2,.,x_{n})$ and $y_0 = (x_2',.,x_{n}')$:
\begin{eqnarray}
|b^2(y_0) - b^2(y)| = |b^2(x_2',.,x_{n}') - b^2(x_2,.,x_{n})| &\leq& \displaystyle{\sum_{i=3}^{n} |b^2(x_2',.,x_{i}',x_{i+1},.,x_{n}) - b^2(x_2',.,x_{i-1}',x_{i},.,x_{n})|} \nonumber \\
&\leq& \displaystyle{\sum_{i=3}^{n} |\frac{\partial (b^2)}{\partial x_i}|^{D^t}_{0} |x_{i}' - x_{i}|} \leq (n-2) b^* |y' - y|  \nonumber 
\end{eqnarray}
Now, let $(t,S)$ and $(t',S')$ belong to {$]0,T[$ \texttimes \hspace{.02cm} $\Omega_S$}. We compute
\begin{eqnarray} 
|\int_{\Omega_y^S}{p(t,S,y) (b^2(y_0) - b^2(y)) dy} - \int_{\Omega_y^{S'}}{p(t',S',y) (b^2(y_0) - b^2(y)) dy}|  \nonumber \\
\leq \int_{\Omega_y^S \cap \Omega_y^S}{|p(t,S,y) - p(t',S',y)| |b^2(y_0) - b^2(y)| dy} + \int_{\Omega_y^S \setminus \Omega_y^{S'}}{p(t,S,y) |b^2(y_0) - b^2(y)| dy}  \nonumber \\
+ \int_{\Omega_y^{S'} \setminus \Omega_y^S}{p(t,S,y) |b^2(y_0) - b^2(y)| dy} \nonumber \\
\leq H^{D^T}_{h} (p) D((t,S),(t',S')) \int_{\Omega_y^S \cap \Omega_y^S}{|b^2(y_0) - b^2(y)| dy} + |p|^{D^T}_{0} \int_{\Omega_y^S \setminus \Omega_y^{S'}}{|b^2(y_0) - b^2(y)| dy} \nonumber \\
+ |p|^{D^T}_{0} \int_{\Omega_y^{S'} \setminus \Omega_y^S}{|b^2(y_0) - b^2(y)| dy} \nonumber \\
\leq b^* K |p|^{D^T}_{h} (D((t,S),(t',S')) + \int_{\Omega_y^S \setminus \Omega_y^{S'}}{|y_0 - y| dy} + \int_{\Omega_y^{S'} \setminus \Omega_y^S}{|y_0 - y| dy}) \nonumber
\end{eqnarray}
By assumption on the boundary of our domain, it is possible to find a constant K depending only on $\Omega$ such as $\forall S,S' \in \Omega_S$, $\int_{\Omega_y^S \setminus \Omega_y^{S'}}{|y_0 - y| dy} \leq K D(S,S')$. This gives us
\begin{eqnarray}
|\int{p(\underline {b}^2 - b^2)dy}|^{D^T}_{h} \leq b^* K |p|^{D^t}_{h} \nonumber
\end{eqnarray}
And since $|\int{p(\underline {b}^2 - b^2)dy}|^{D^T}_{2+h} = |\int{p(\underline {b}^2 - b^2)dy}|^{D^T}_{h} + |\int{\frac{\partial p}{\partial t}(\underline {b}^2 - b^2)dy}|^{D^T}_{h} + |\int{\frac{\partial p}{\partial S}(\underline {b}^2 - b^2)dy}|^{D^T}_{h} + |\int{\frac{\partial^2 p}{\partial S^2}(\underline {b}^2 - b^2)dy}|^{D^T}_{h}$, we get from the previous computation
\begin{eqnarray}
|\int{p(\underline {b}^2 - b^2)dy}|^{D^T}_{2+h} \leq b^* K |p|^{D^t}_{2+h} \nonumber
\end{eqnarray}

\noindent
We now have to find a bound on $|\frac{1}{\int{b^2pdy}}|^{D^T}_{2+h}$. Since $p$ belongs to $X^{T^*}_{x^*}$, we have {$|\frac{1}{\int{b^2pdy}}|^{D^T}_{0} \leq \frac{2}{\delta_1^2 \underline{p_0} V(\Omega)}$}. Now, let $(t,S)$ and $(t',S')$ belong to {$]0,T[$ \texttimes \hspace{.02cm} $\Omega_S$}. We write
\begin{eqnarray}
|\frac{1}{\int_{\Omega_y^S}{b^2(y) p(S,t,y)dy}} - \frac{1}{\int_{\Omega_y^{S'}}{b^2(y) p(S',t',y)dy}}| &\leq&
\frac{|\int_{\Omega_y^{S'}}{b^2(y) p(S',t',y)dy} - \int_{\Omega_y^S}{b^2(y) p(S,t,y)dy}|}{\int_{\Omega_y^S}{b^2(y) p(S,t,y)dy} \int_{\Omega_y^{S'}}{b^2(y) p(S',t',y)dy}} \nonumber \\
&\leq& K H^{D^t}_{h} (p) D((t,S),(t',S')) \nonumber
\end{eqnarray}
We used the same kind of arguments than earlier, K denotes here another constant depending on $\delta_1$, $\delta_2$, $\underline{p_0}$ and $\Omega$. This gives us $|\frac{1}{\int{b^2pdy}}|^{D^T}_{h} \leq K (1 + |p|^{D^T}_{h})$. 
As for derivatives of $\frac{1}{\int{b^2pdy}}$, for instance with respect to $S$, we have
\begin{eqnarray}
|\frac{\partial}{\partial S} (\frac{1}{\int{b^2pdy}})|^{D^T}_{h} = |-\frac{\int{b^2 \frac{\partial p}{\partial S}dy}}{(\int{b^2pdy})^2}|^{D^T}_{h} \leq K (1 + |p|^{D^T}_{h})^2 |p|^{D^T}_{1+h} \nonumber
\end{eqnarray}
The same kind of computation is true for the derivative of second order 
\begin{eqnarray}
|\frac{\partial}{\partial S} (-\frac{\int{b^2 \frac{\partial p}{\partial S}dy}}{(\int{b^2pdy})^2})|^{D^T}_{h} = |\frac{2 (\int{b^2 \frac{\partial p}{\partial S}dy})^2} {(\int{b^2pdy})^3} - \frac{\int{b^2 \frac{\partial^2 p}{\partial S^2}dy}} {(\int{b^2pdy})^2}|^{D^T}_{h} \leq K [(1 + |p|^{D^T}_{h})^3 (|p|^{D^T}_{1+h})^2 + (1 + |p|^{D^T}_{h})^2 |p|^{D^T}_{2+h}] \nonumber
\end{eqnarray}

\noindent
Eventually, we get
\begin{eqnarray}
|\frac{1}{\int{b^2pdy}}|^{D^T}_{2+h} \leq K (1 + |p|^{D^T}_{2+h} + (|p|^{D^T}_{2+h})^2 + (|p|^{D^T}_{2+h})^3 + (|p|^{D^T}_{2+h})^4 + (|p|^{D^T}_{2+h})^5) \nonumber
\end{eqnarray}

\noindent
Combining this result with the previous computations, we find 
\begin{eqnarray}
|I(p) - \frac{1}{{\underline {b}}^2}|^{D^T}_{2+h} \leq b^* K P(|p|^{D^t}_{2+h}) \nonumber
\end{eqnarray} 
with P a polynomial function of degree 6, strictly positive on $\mathbb{R}_+^*$.
Now, writing $\sqrt{I(p)} - \frac{1}{\underline {b}} =  \frac{I(p) - \frac{1}{{\underline {b}}^2}}{\sqrt{I(p)} + \frac{1}{\underline {b}}}$, we find the same kind of results for the second term involved in the lemma. This concludes the proof.

\section{Conclusion}
In this paper, we have shown that the equation driving the calibration problem for local and stochastic volatility models is well-posed in the case of suitably regularized initial conditions. It is however clear that the solution of the full Kolmogorov equation with Dirac initial condition does not obtain as a consequence of Theorem 1 : possible extensions of our results towards this direction are currently being explored. Let us also mention that a generalization of Theorem 1 to multidimensional correlation calibration have already been investigated and will be presented in \cite{tachet}

\end{document}